\documentclass[reqno,fleqn,onecolumn]{amsproc}%
\usepackage{amsfonts}
\usepackage{amsmath}
\usepackage{amssymb}
\usepackage{graphicx}
\usepackage[none]{hyphenat}%
\theoremstyle{plain}

\numberwithin{equation}{section}
\begin{document}
\title{}

\begin{center}
\textbf{Translation Invariance of Investment}

\bigskip

Yukio Hirashita

\bigskip

\textit{Faculty of Liberal Arts, Chukyo University}

\textit{ Nagoya, 466-8666, Japan }

\textit{yukioh@cnc.chukyo-u.ac.jp}

\bigskip
\end{center}

\begin{quotation}
For a game with positive profit, the optimal proportion of investment required
to continue investing without borrowing is uniquely determined by an integral
equation for each price. For a game with parallel translated profit, the ratio
of the optimal proportion of investment to its price has some invariance
properties. The optimal price of a game with parallel translated profit
converges to its expectation divided by $e$ to the riskless interest rate for
a certain period.\bigskip



\bigskip

\bigskip
\end{quotation}

\noindent\textbf{1. Introduction}

\noindent A game $(a(x),$ $F(x))$ would imply that if the investor invests $1$
dollar, then he or she receives $a(x)$ dollars (including the $1$ dollar
initially invested) in accordance with a distribution function $F(x)$ that is
defined on an interval $I\subseteq(-\infty$, $\infty)$ such that $\int
_{I}d(F(x))=1$. For each positive price, in order to continue investing
without borrowing, the optimal proportion of investment must be in $(0$, $1]$
(see [4,5,6,7]). In Section 2, we will investigate a game with parallel
translated profit $(a(x)+n,$ $F(x))$.

For simplicity, we will omit the currency notation. It is assumed that the
profit function $a(x)$ is measurable and non constant (a.e.) with respect to
$F(x)$. In the absence of confusion, we use $dF$ to represent $d(F(x))$ and
employ the following notation: $E:=\int_{I}a(x)dF,$ $H:=\int_{I}1/a(x)dF,$
$\xi:=\mathrm{ess}$ $\inf_{x\in I}\,a(x),$ $H_{\xi}:=\int_{I}1/(a(x)-\xi)dF$.
In this paper, we always assume that $\xi>0$. If $\int_{a(x)=\xi}dF>0$, we
define $H_{\xi}$ $=\infty$ and $1/H_{\xi}$ $=0.$ Since $a(x)$ is non-constant,
we have $\xi\leq\xi+1/H_{\xi}<1/H<E$ (see [4, Section 2 and Lemma 3.5]).

For each $u\in(\xi+1/H_{\xi}$, $E)$, the pre-optimal proportion $\widetilde
{t}_{u}\in(0$, $u/(u-\xi))$ is uniquely determined by the equation $\int
_{I}(a(x)-u)/(a(x)\widetilde{t}_{u}-u\widetilde{t}_{u}+u)dF=0$ (see [4,
Section 3]). $\widetilde{t}_{u}$ has the following properties: (1)
$\widetilde{t}_{u}$ is strictly decreasing ([4, Lemma 3.6]); (2)
$\lim_{u\rightarrow E^{-}}\widetilde{t}_{u}=0$ ([4, Lemma 3.7]); (3)
$\widetilde{t}_{1/H}=1$ ([4, Lemma 3.8]); and (4) $\lim_{u\rightarrow
(\xi+1/H_{\xi})^{+}}\widetilde{t}_{u}=1+\xi H_{\xi}$ ([4, Lemma 3.9]).

For each $u\in(\xi+1/H_{\xi}$, $E)$ and $t\in(0$, $u/(u-\xi))$, the pre-growth
rate $\widetilde{G}_{u}(t)$ is defined by $\exp\left(  \int_{I}\log\left(
a(x)t/u-t+1\right)  dF)\right)  $ (see [4, Section 4 ]). We always assume that
$\widetilde{G}<\infty$, which implies that $\int_{a(x)>1}\log a(x)dF<+\infty$
(see [4, Lemma 4.3]). $\widetilde{G}_{u}(\widetilde{t}_{u})$ has the following
properties: (1) $\widetilde{G}_{u}(\widetilde{t}_{u})$ is continuous and
strictly decreasing ([4, Lemma 4.9 and Theorem 4.1]); (2) $\lim_{u\rightarrow
E^{-}}\widetilde{G}_{u}(\widetilde{t}_{u})=1$ ([4, Lemma 4.13]); (3)
$\widetilde{G}_{1/H}(\widetilde{t}_{1/H})$ $=H\exp\left(  \int_{I}\log
a(x)dF\right)  $ ([4, Lemma 4.21]); and (4) $\lim_{u\rightarrow(\xi+1/H_{\xi
})^{+}}\widetilde{G}_{u}(\widetilde{t}_{u})$ $=H_{\xi}\exp(\int_{I}\log\left(
a(x)-\xi\right)  dF)$ ([4, Lemmas 3.9 and 4.20]).

We cite the following lemmas (for the proofs see [4, Section 5]):\bigskip

\noindent\textbf{Lemma 1.1 }[4, Corollary 5.1]\textbf{.} \textit{For each
}$u\in(1/H,$\textit{ }$E)$\textit{, the optimal proportion of investment is
equal to }$\widetilde{t}_{u}$\textit{, and the maximized limit expectation of
the growth rate is equal to }$\widetilde{G}_{u}(\widetilde{t}_{u})$%
\textit{.}\bigskip

\noindent\textbf{Lemma 1.2 }[4, Corollary 5.3]\textbf{.} \textit{For each
}$u\in(0,$\textit{ }$1/H$\textit{ }$]$\textit{, the optimal proportion of
investment is equal to }$1$\textit{, and the maximized limit expectation of
the growth rate is equal to }$\exp(\int_{I}\log a(x)$\textit{ }$dF)/u$%
\textit{.}\bigskip

If the riskless continuously compounded interest rate $r>0$ for a period is
given, the optimal price of a game can be uniquely determined by the equation
such that the maximized limit expectation of the growth rate is equal to
$e^{r}$ (see [4, Section 6]).

\bigskip

\bigskip

\noindent2\textbf{. Parallel translated profit}

\noindent For a game with parallel translated profit $(a(x)+n,$ $F(x))$
$(n>-\xi)$, we use underlined notations such as $\underline{a}(x)=a(x)+n,$
$\underline{E}=E+n,$ $\underline{H}=\int_{I}1/\underline{a}(x)dF,$
$\underline{\xi}=\xi+n$ and $\underline{H}_{\underline{\xi}}=\int_{I}1/\left(
\underline{a}(x)-\underline{\xi}\right)  dF$.\bigskip

\noindent\textbf{Example 1}. Suppose that the profit $a$ or $b$ ($0<b<a$)
occurs with the probability $p$ or $1-p$ $(0<p<1)$, respectively. From%

\begin{align}
&  \int_{I}\frac{(a(x)+n)-(u+n)}{(a(x)+n)\widetilde{\underline{t}}%
_{u+n}-(u+n)\widetilde{\underline{t}}_{u+n}+u+n}dF\nonumber\\
&  =\frac{a-u}{(a-u)\widetilde{\underline{t}}_{u+n}+u+n}p+\frac{b-u}%
{(b-u)\widetilde{\underline{t}}_{u+n}+u+n}(1-p)=0, \tag{2.1}%
\end{align}
we have $\widetilde{\underline{t}}_{u+n}=(E-u)(n+u)/((a-u)(u-b))$ $(b<u<E<a$,
$-b<n)$. Thus, we obtain%
\begin{align}
\underline{\widetilde{G}}_{u+n}(\widetilde{\underline{t}}_{u+n})  &
=\exp\left(  \int_{I}\log\left(  \frac{(a(x)+n)\widetilde{\underline{t}}%
_{u+n}}{u+n}-\widetilde{\underline{t}}_{u+n}+1\right)  dF\right) \nonumber\\
&  =(a-b)\left(  \frac{p}{u-b}\right)  ^{p}\left(  \frac{1-p}{a-u}\right)
^{1-p}, \tag{2.2}%
\end{align}
which is independent of $n$.\bigskip

\noindent\textbf{Lemma 2.1.} $\widetilde{G}<+\infty$, \textit{if and only if}
$\widetilde{\underline{G}}<+\infty$.

\noindent\textbf{\textbf{Proof}. }If $\widetilde{G}<+\infty$, then
from\textbf{ }$\int_{a(x)+n>1}\log(a(x)+n)dF=\int_{1-n<a(x)\leq|n|+1}%
\log(a(x)+n)dF$ $+\int_{a(x)>|n|+1}\log(a(x)+n)dF\leq$ $\log(2|n|+1)+\log
2+\int_{a(x)>1}\log a(x)dF$, we have $\widetilde{\underline{G}}<+\infty$. The
remaining can be verified in a similar manner.\hfill$\square$\bigskip

\noindent\textbf{Theorem 2.2.} $\widetilde{\underline{t}}_{u+n}%
/(u+n)=\widetilde{\underline{t}}_{u}/u$\textit{ }$(u\in(\xi+1/H_{\xi}%
$\textit{, }$E))$\textit{ is independent of }$n>-\xi$.

\noindent\textbf{Proof.} By definition, we have $\int_{I}\left(
a(x)-u\right)  /\left(  (a(x)-u)\widetilde{\underline{t}}_{u+n}+u+n\right)
dF=0$, which is equivalent to $\int_{I}\left(  u+n\right)  /\left(
(a(x)-u)\widetilde{\underline{t}}_{u+n}+u+n\right)  dF$ $=1$. By
differentiating them, we obtain%
\begin{equation}
\frac{\partial\widetilde{\underline{t}}_{u+n}}{\partial n}=-\frac{\int
_{I}\frac{a(x)-u}{\left(  (a(x)-u)\widetilde{\underline{t}}_{u+n}+u+n\right)
^{2}}dF}{\int_{I}\frac{\left(  a(x)-u\right)  ^{2}}{\left(  (a(x)-u)\widetilde
{\underline{t}}_{u+n}+u+n\right)  ^{2}}dF}, \tag{2.3}%
\end{equation}%
\begin{equation}
(\widetilde{\underline{t}}_{u+n}-(u+n)\frac{\partial\widetilde{\underline{t}%
}_{u+n}}{\partial n})\int_{I}\frac{a(x)-u}{\left(  (a(x)-u)\widetilde
{\underline{t}}_{u+n}+u+n\right)  ^{2}}dF=0. \tag{2.4}%
\end{equation}
Since these functions are analytic with respect to $u$ (see [4, Lemma 3.1] and
[8, Chapter VIII]), we have $\partial\widetilde{\underline{t}}_{u+n}/\partial
n$ $-\widetilde{\underline{t}}_{u+n}/(u+n)=0$ or $\partial\widetilde
{\underline{t}}_{u+n}/\partial n=0$ as a constant function. If $\partial
\widetilde{\underline{t}}_{u+n}/\partial n$ $=\widetilde{\underline{t}}%
_{u+n}/(u+n)$, then $\partial\left(  \widetilde{\underline{t}}_{u+n}%
/(u+n)\right)  /\partial n$ $=\left(  (u+n)\partial\widetilde{\underline{t}%
}_{u+n}/\partial n-\widetilde{\underline{t}}_{u+n}\right)  /(u+n)^{2}=0$,
which implies $\widetilde{\underline{t}}_{u+n}/(u+n)=\widetilde{\underline{t}%
}_{u}/u$. If $\partial\widetilde{\underline{t}}_{u+n}/\partial n=0$, then
$\widetilde{\underline{t}}_{u+n}=\widetilde{t}_{u}$, particularly,
$\widetilde{\underline{t}}_{1/H+n}=\widetilde{t}_{1/H}$. Since $\widetilde
{\underline{t}}_{1/H+n}=\widetilde{t}_{1/H}=1=\underline{\widetilde{t}%
}_{1/\underline{H}}$ and $\widetilde{\underline{t}}_{u+n}$ is strictly
decreasing with respect to $u$, we have $1/\underline{H}=1/H+n$, which
contradicts the fact that $1/\underline{H}-n$ is strictly increasing with
respect to $n$ (see [5, the proof of Lemma 1.1]).\hfill$\square$\bigskip

\noindent\textbf{Theorem 2.3.}\textit{ }$\underline{\widetilde{G}}%
_{u+n}(\widetilde{\underline{t}}_{u+n})=\widetilde{G}_{u}(\widetilde{t}_{u}%
)$\textit{ }$(u\in(\xi+1/H_{\xi}$\textit{, }$E))$\textit{ is independent of
}$n>-\xi$\textit{.}

\noindent\textbf{Proof.} By differentiating $\underline{\widetilde{G}}%
_{u+n}(\widetilde{\underline{t}}_{u+n})$, we have
\begin{equation}
\frac{\partial}{\partial n}\underline{\widetilde{G}}_{u+n}(\widetilde
{\underline{t}}_{u+n})=\left(  \frac{\partial\widetilde{\underline{t}}_{u+n}%
}{\partial n}-\frac{\widetilde{\underline{t}}_{u+n}}{u+n}\right)
\underline{\widetilde{G}}_{u+n}(\widetilde{\underline{t}}_{u+n})\int_{I}%
\frac{a(x)-u}{(a(x)-u)\widetilde{\underline{t}}_{u+n}+u+n}dF. \tag{2.5}%
\end{equation}
Therefore, the conclusion is obtained form the proof of Theorem 2.2.\hfill
$\square$\bigskip

It should be noted that the difference $\underline{E}-(\underline{\xi
}+1/\underline{H}_{\underline{\xi}})=E-\left(  \xi+1/H_{\xi}\right)  $ is
independent of $n$.\bigskip

\noindent\textbf{Lemma 2.4.} \textit{If }$E<\infty$\textit{, }$\lim
_{n\rightarrow\infty}\left(  \underline{E}-1/\underline{H}\right)
=0$\textit{.}

\noindent\textbf{\textbf{Proof}.} By definition, we have%
\begin{align}
\underline{E}-1/\underline{H}  &  =E+n-\frac{1}{\int_{I}\frac{1}{a(x)+n}%
dF}\nonumber\\
&  =E-\left(  \frac{1}{\int_{I}\frac{1}{a(x)+n}dF}-n\right)  \;(n>-\xi).
\tag{2.6}%
\end{align}
Thus, using [5, Lemma 1.1], we derive the conclusion.\hfill$\square$\bigskip

\noindent\textbf{Lemma 2.5.} $\lim_{n\rightarrow\infty}\underline
{\widetilde{G}}_{1/\underline{H}}(\underline{\widetilde{t}}_{1/\underline{H}%
})=1$.

\noindent\textbf{\textbf{Proof}.} The conclusion is derived from the fact that
$\underline{\widetilde{G}}_{1/\underline{H}}(\underline{\widetilde{t}%
}_{1/\underline{H}})=\underline{H}\exp\left(  \int_{I}\log\left(
a(x)+n\right)  dF\right)  $, $\lim_{n\rightarrow\infty}n\underline{H}%
=\lim_{n\rightarrow\infty}\int_{I}n/(a(x)+n)dF=1$, and
\begin{equation}
\lim_{n\rightarrow\infty}\frac{\exp\left(  \int_{I}\log\left(  a(x)+n\right)
dF\right)  }{n}=\lim_{n\rightarrow\infty}\exp\left(  \int_{I}\log\left(
\frac{a(x)}{n}+1\right)  dF\right)  =1. \tag{2.7}%
\end{equation}
\hfill$\square$\bigskip

Let $r>0$ be the riskless continuously compounded interest rate such that
$e^{r}<\widetilde{G}_{1/H}(1)$. The optimal price $u_{r}\in(1/H$, $E)$ of a
game $(a(x),$ $F(x))$ is then uniquely determined by the equation
$e^{r}=\widetilde{G}_{u_{r}}(\widetilde{t}_{u_{r}})$ (see Lemma 1.1 and [4,
Section 6]).\bigskip

\noindent\textbf{Theorem 2.6. }$\underline{u}_{r}=u_{r}+n$\textit{, if }%
$e^{r}<\min(\widetilde{G}_{1/H}(1)$\textit{, }$\underline{\widetilde{G}%
}_{1/\underline{H}}(1))$\textit{ and }$n>-\xi$.

\noindent\textbf{\textbf{Proof}. }By definition we have $\widetilde{G}_{u_{r}%
}(\widetilde{t}_{u_{r}})=\underline{\widetilde{G}}_{\underline{u}_{r}%
}(\widetilde{\underline{t}}_{\underline{u}_{r}})=e^{r}$. From Theorem 2.3, we
observe that $\underline{\widetilde{G}}_{u_{r}+n}(\widetilde{\underline{t}%
}_{u_{r}+n})$ $=\widetilde{G}_{u_{r}}(\widetilde{t}_{u_{r}})$. Based on the
uniqueness of price $\underline{u}_{r}$, we obtain $\underline{u}_{r}=u_{r}%
+n$.\hfill$\square$\bigskip

\noindent\textbf{Theorem 2.7.}\textit{ If }$E<\infty$\textit{, }%
$\lim_{n\rightarrow\infty}\underline{u}_{r}/\underline{E}=1/e^{r}$\textit{.}

\noindent\textbf{\textbf{Proof}.} From $\lim_{n\rightarrow\infty}%
\underline{\widetilde{G}}_{1/\underline{H}}(\underline{\widetilde{t}%
}_{1/\underline{H}})=1$ (Lemma 2.5), for each sufficiently large $n$, we have
$e^{r}>\underline{\widetilde{G}}_{1/\underline{H}}(1)$. In this case, the
price $\underline{u}_{r}<1/\underline{H}$ is uniquely determined by the
equation $e^{r}=\exp(\int_{I}\log\left(  a(x)+n\right)  dF)/\underline{u}_{r}$
(see Lemma 1.2 and [4, Section 6]). Thus, we obtain%
\[
\lim_{n\rightarrow\infty}\frac{\underline{u}_{r}}{\underline{E}}%
=\lim_{n\rightarrow\infty}\left(  \exp\left(  \int_{I}\log\left(  \frac
{a(x)}{n}+1\right)  dF\right)  \times\frac{n}{E+n}\right)  \frac{1}{e^{r}%
}=\frac{1}{e^{r}}.\;\square
\]
\bigskip

\noindent\textbf{Example 2}. Suppose that the profit $1$ or $19$ occurs with
each probability $1/2$. Then, using Example 1, we have $E=10,$ $H=10/19,$
$\xi=1,$ $H_{\xi}=\infty$, $\underline{E}=10+n,$ $\underline{H}=(n+10)/((n+1)$
$(n+19)),$ $\underline{\widetilde{G}}_{u+n}(\widetilde{\underline{t}}%
_{u+n})=9/\sqrt{(u-1)(19-u)}$, $\underline{\widetilde{G}}_{1/\underline{H}%
}(\underline{\widetilde{t}}_{1/\underline{H}})=(n+10)/\sqrt{(n+1)(n+19)}$,
$n>-1$, and $1<u<10.$

Let the riskless continuously compounded interest rate $r$ be $0.05$. Then, as
$e^{r}\fallingdotseq1.051<\widetilde{G}_{1/H}(1)\fallingdotseq2.294$, the
equation $e^{r}=9/\sqrt{(u-1)(19-u)}$ deduces the optimal price $u_{r}%
=10-9\sqrt{1-e^{-2r}}\fallingdotseq7.224$, which is considerably different
from $E/e^{r}\fallingdotseq9.512$.

By the equality $e^{r}=(n_{0}+10)/\sqrt{(n_{0}+1)(n_{0}+19)}$, we have
$n_{0}=9e^{r}/\sqrt{e^{2r}-1}$ $-10\fallingdotseq19.175$. Therefore, if
$-1<n<19.175$, then $\underline{u}_{r}\fallingdotseq7.224+n$ (Theorem 2.6).

In the case where $n=99>19.175$, the optimal price $\underline{u}_{r}$ is
determined by $e^{r}=\exp$ $(\int_{I}\log\left(  a(x)+n\right)  dF)/\underline
{u}_{r}$. Thus, $\underline{u}_{r}=\sqrt{(n+1)(n+19)}/e^{r}\fallingdotseq
103.330,$ which is considerably close to $\underline{E}/e^{r}\fallingdotseq
103.684$, as suggested by Theorem 2.7.\bigskip

\bigskip

\bigskip

\noindent\textbf{References\smallskip}

\noindent\lbrack1] N. H. Bingham, C. M. Goldie, J. L. Teugels, \textit{Regular
variation}, Cambridge

University Press, Cambridge, 1987.

\noindent\lbrack2] F. Black, M. Scholes. \textit{The pricing of options and
corporate liabilities}.

Journal of Political Economy 81(1973), 637-54.

\noindent\lbrack3] W. Feller, \textit{An introduction to probability theory
and \ its application}, vol. I, II,

John Wiley and Sons, New York, 1957, 1966.

\noindent\lbrack4] Y. Hirashita, \textit{Game pricing and double sequence of
random variables},

Preprint, arXiv:math.OC/0703076 (2007).

\noindent\lbrack5] Y. Hirashita, \textit{Ratio of price to expectation and
complete Bernstein function},

Preprint, arXiv:math.OC/0703077 (2007).

\noindent\lbrack6] J. L. Kelly, \textit{A new interpretation of information
rate},\ Bell System Tech. J.

35 (1956) 917--926.

\noindent\lbrack7] D. G. Luenberger, \textit{Investment science}, Oxford
University Press, Oxford, 1998.

\noindent\lbrack8] D. V. Widder, \textit{The Laplace Transform}, Princeton
University Press, Princeton,

1941.
\end{document}